# 4D Dyck triangle and its projections

Gennady Eremin
ergenns@gmail.com

July 21, 2020

**Abstract**. The classic Dyck triangle, the Catalan triangle, and the Catalan convolution matrix are plane projections of the multidimensional Dyck triangle. In the Dyck path, each node is uniquely determined by two of four interrelated parameters: (*i*) the position of the current parenthesis, (*ii*) the current unbalance of the parentheses, (*iii*) the number of viewed left parentheses, and (*iv*) the same for right parentheses. The last two parameters can be redefined, respectively, as the index of the current Catalan number and the index of the summand in the decomposition of the Catalan number into the sum of squares (Dyck squares). For the 4D Dyck triangle, we consider six 2D projections (some of them are not yet in demand) and four 3D projections.

## 1  Introduction

This paper continues [1]. A Dyck path of length $2n$ (*semilength n*, *size n*) is a diagonal lattice path from the origin to node $(2n, 0)$ consisting of $n$ upsteps in the form of vectors $(1, 1)$ and $n$ downsteps in the form of vectors $(1, -1)$, such that the path never goes below the ground level. In the bracket set, the left parenthesis corresponds to the upstep, and the right parenthesis is the downstep.

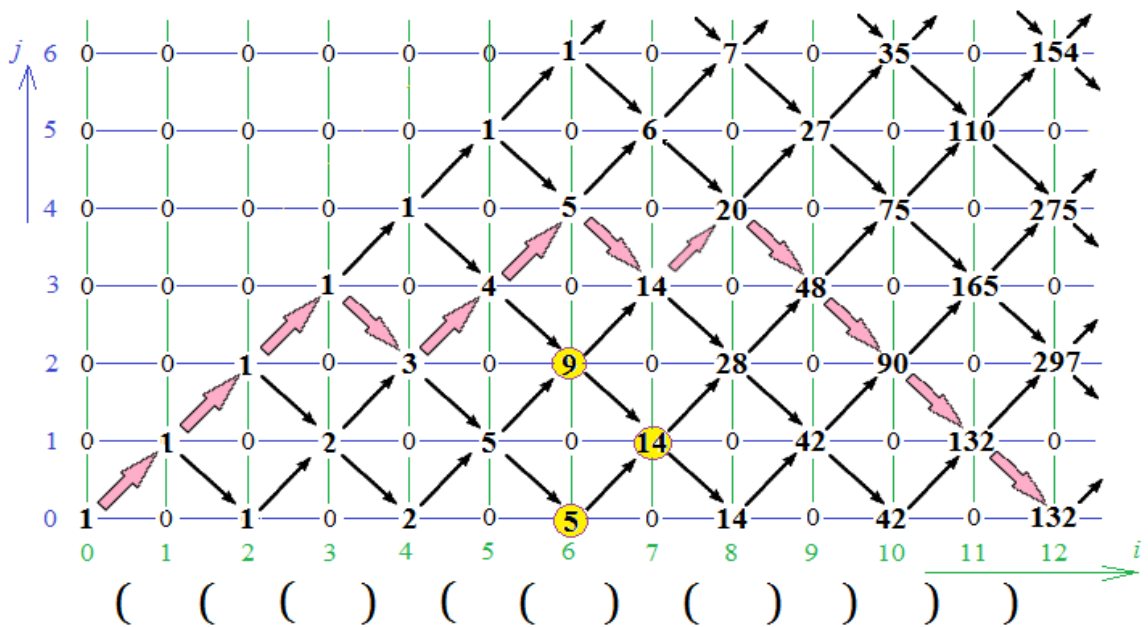

Figure 1: The Dyck path and the corresponding sequence of parentheses.

Figure 1 collected several pictures from [3]. In the grid, the *positions* of the diagonal vectors are plotted along the abscissa axis, *i-axis*, and the *unbalance* (exceeding the number of upsteps above downsteps) is along the ordinate axis, *j-axis*. Arrows show valid vectors. Diagonal vectors pass accessible nodes, for which the sum of coordinates is even. In the *i*th position, the unbalance $j$ cannot exceed $i$, so we have a triangular design in the $i \times j$ grid, which we will call a *Dyck ij-triangle*. Thus, the node $(i, j)$, *ij-node*, belongs to the Dyck *ij*-triangle if and only if

$$j \leq i \text{ and } i+j = \text{even}. \tag{1}$$

Each *ij-node* has a label $d(i, j)$, *dynamics*, which is equal to the number of paths starting at $(0, 0)$ and ending at $(i, j)$. For any *ij*-node, the dynamics is determined as follows:
$$d(0,0) = 1, \ d(i,j) = d(i–1, j+1) + d(i–1, j–1). \tag{2}$$

The dynamics is equal to zero in unreachable nodes (condition (1) is not satisfied). Three nodes connected by (2) are highlighted in yellow: $d(7,1) = d(6,2) + d(6,0)$.

In the lower two lines, we have the Catalan numbers: $d(2n, 0) = d(2n–1, 1) = \text{Cat}(n)$. For clarity, Figure 1 shows the Dyck *ij-path* by thick colored arrows, below the *i*-axis there is a corresponding parentheses (the Dyck word of semilength 6). We will repeat this *ij*-path and three selected *ij*-nodes in the following figures.

## 2  Isolines of Dyck triangle

Any rectangular coordinate system has horizontal and vertical lines, *isolines*. In Figure 1, horizontal lines are *j-isolines*, since all nodes (accessible and inaccessible) of the horizontal have the same $j$ coordinate. Accordingly, the vertical lines are *i-isolines*, since the nodes of any vertical have the same $i$ coordinate. In this section, we will consider additional isolines in the Dyck *ij*-triangle.

Accordance to (1) for a reachable node $(i, j)$ the value of $i+j$ is even. Let's write it like this: $i+j = 2n$. Consider a small example.

**Example 1**. In the Dyck *ij*-path (see Figure 2), the last four downward vectors pass through five nodes on a diagonal section from $(8, 4)$ to $(12, 0)$. For nodes of the falling diagonal starting at the top $(6, 6)$, the sum of the coordinates is $2 \times 6 = 12$. The dynamics of the lower diagonal point is equal to the 6th Catalan number, Cat $(6) = d(12, 0) = 132$, and it is logical to connect such a diagonal with the 6th Catalan number, or rather with the 6th index. For example, we can say that the node $(7, 3)$ is tied to Cat $(5)$, because $(7+3)/2 = 5$.  □

The falling diagonal of the Dyck *ij*-triangle will be called an *n-isoline*. But in the triangle, there are also rising diagonals. It is easy to see that condition (1) implies the following: for any reachable node $(i, j)$ the value of $i−j$ is even or zero. Let's write it this way: $i−j = 2k$. Consider another example.

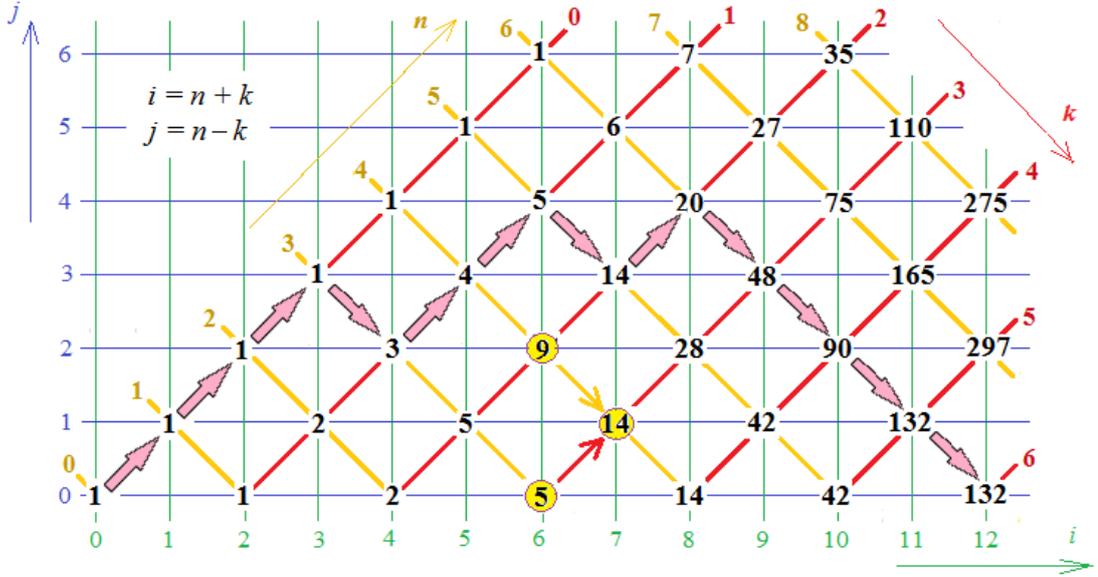

Figure 2: The Dyck ij-triangle.

**Example 2**. In the same Dyck *ij*-path, the first three upsteps pass through four nodes from (0, 0) to (3, 3) at the beginning of the main diagonal. For nodes of this rising diagonal, $k = 0$. The following ascending vectors pass through points (4, 2), (5, 3), and (6, 4) for which $k = 1$. In other words, we can say that these points are on the rising diagonal #1. For an achievable node, the ordinal number of the corresponding rising diagonal is equal to the half-difference of the coordinates. □

The rising diagonal of the Dyck *ij*-triangle will be called a *k-isoline*. Obviously, any accessible node $(v, w)$ is located at the intersection of four isolines: (*i*) the *v*-th *i*-isoline (vertical), (*ii*) the *w*-th *j*-isoline (horizontal), (*iii*) the $(v + w)/2$-th *n*-isoline (falling diagonal), and (*iv*) the $(v - w)/2$-th *k*-isoline (rising diagonal). Four 0-th isolines intersect at the origin. Thus, for an attainable node $(i, j)$

$$i+j = 2n, \quad i-j = 2k \quad \text{or} \quad i = n+k, \; j = n-k \quad (i, j, n, k \geq 0). \tag{3}$$

In Figure 2, unreachable nodes are removed, and *n*-isolines from #0 to #8 are shown in yellow. The yellow arrow shows the direction of the virtual *n*-axis. In red, we showed the *k*-isolines from #0 to #6. The red arrow shows the direction of the virtual *k*-axis. It is easy to see that inequalities $i \geq n \geq j$, $n \geq k$ follow from (3). This is easy to check from the picture. Here we have repeated the previous Dyck *ij*-path and the three selected *ij*-nodes. We pay attention to the fact that the Dyck paths pass only diagonals, that is, *n*-isolines and *k*-isolines. The upsteps lie on the *k*-isolines, and the downsteps lie on the *n*-isolines.

Let's give a concrete meaning to the variable $k$ from (3). Any node $(v, v)$, $v > 0$, is the intersection of *i*-isoline #$v$, *j*-isoline #$v$, *n*-isoline #$v$, and *k*-isoline #0. It is not easy to obtain the following equality:

$$\text{Cat}(v) = d^2(v, v) + d^2(v, v-2) + d^2(v, v-4) + \ldots + d^2(v, v-2\lfloor v/2 \rfloor) \tag{4}$$

or

$$\text{Cat}(v) = \sum_{k=0}^{\lfloor v/2 \rfloor} t_k^2(v), \quad t_k(v) = d(v, v-2k) \qquad (4a)$$

In (4a), the terms are logically called *Dyck squares* (see formulas 2 and 7 in [2], as well as the equivalent formula 12 in [8]). Some terms are obvious, for example,

$t_0(v) = d(v,v) = 1; \quad t_1(v) = d(v, v-2) = v-1; \quad t_2(v) = d(v, v-4) = v(v-3)/2;$ (5)
$t_{\lfloor v/2 \rfloor}(v) = d(v, v-2\lfloor v/2 \rfloor) = \text{Cat}(\lceil v/2 \rceil).$

In this case, the variable $k$ is the index of the Dyck square. Thus, all four considered isolines are related to certain attributes (states) of balanced parentheses. Let's call $i, j, n, k$ *coordinate variables*, and hence the dependences (3) are the *coordinate equations*.

So, we are dealing with 4-dimensional object, the *Dyck ijnk-triangle*. As a result, we can write the dynamics equation (2) in a general form (for 4D space):

$$D(0,0,0,0) = 1, \quad D(i,j,n,k) = D(i-1, j+1, n, k-1) + D(i-1, j-1, n-1, k). \qquad (6)$$

The 4D Dyck triangle can have up to six different 2D projections, one of which is presented in Figure 2. In the Dyck *ij*-triangle, the *i*-isoline #0 and the *j*-isoline #0 are coordinate axes, and these axes can be replaced. In the next section, we consider a known matrix, which is another 2D projection of the Dyck *ijnk*-triangle.

## 3   Catalan convolution matrix and Dyck *nj*-triangle

In the Dyck *ij*-triangle, unreachable nodes are ballast. After removing nodes with zero dynamics, the data array is reduced by half. Let's rotate all *n*-isolines to the vertical position around the top points; as a result (see Figure 3) we get the *Dyck nj-triangle*, the projection in the $n \times j$ grid (*i*-axis is replaced by *n*-axis).

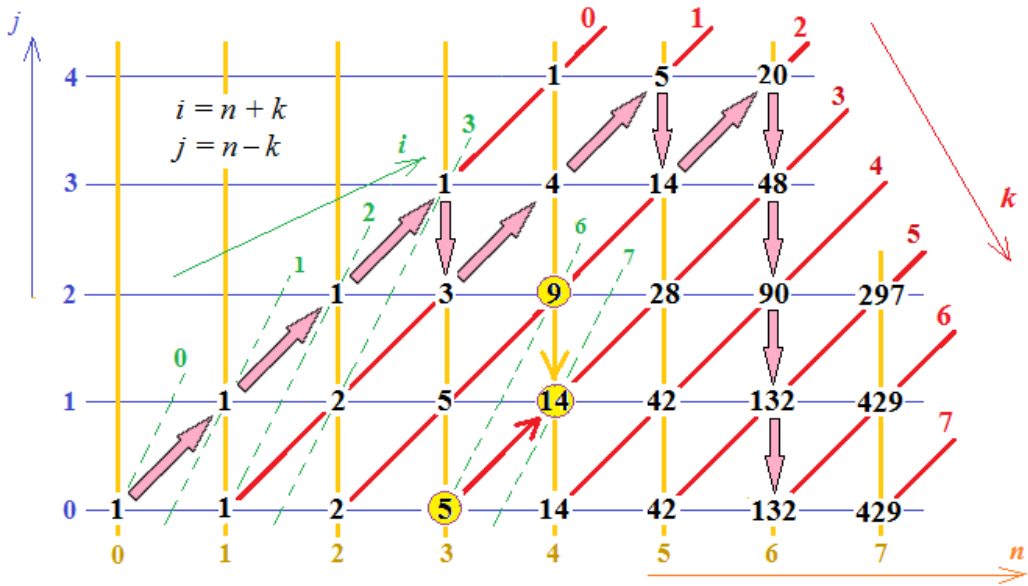

Figure 3: The Dyck *nj*-triangle.

We see that *i*-isolines practically disappeared. For orientation, we have combined some *nj-nodes* by dotted *i*-isolines and showed conditionally the direction of the virtual *i*-axis (recall the coordinate equation $i = n + k$). In the new *nj-path*, the downsteps are transformed into breaks, so the path is reduced by half.

In the new projection, on the abscissa axis we indicate the indexes of the Catalan numbers, so you can not apply (2) to the *nj*-nodes. However, the generalized equation (6) is applicable. For the selected *nj*-nodes, we can write

$$D(7, 1, 4, 3) = D(6, 2, 4, 2) + D(6, 0, 3, 3). \tag{7}$$

The equality (7) is universal for all projections, including the 4D Dyck triangle.

Let's define the relationship between the items of the Dyck *ij*-triangle and the Dyck *nj*-triangle (or between *ij*-nodes and *nj*-nodes). The Dyck *nj*-triangle is known in the literature as the *Catalan convolution matrix* [4, 5]. Elements of this matrix (and correspondingly elements of the Dyck *nj*-triangle) are defined as follows (see [5], p. 2928):

$$C(n, j) = \binom{2n - j}{n - j} - \binom{2n - j}{n - j - 1}, \quad n, \ j \geq 0.$$

Next, let us use the coordinate equations (3) to obtain the general term of the Dyck square (in addition to (5)):

$$t_k(i) = C(n, j) = C(i - k, i - 2k) = \binom{i}{k} - \binom{i}{k - 1} \tag{8}$$

A detailed analysis of (8) is made in [2]. As we see, the use of different projections of the 4D Dyck triangle, can simplify some calculations.

## 4  Lattice paths in $n \times n$ grid and Dyck *nk*-triangle

In the considered projections (see Figure 2 and Figure 3), at each *ijnk*-node of the form $(i, j = 0, n, k)$ the last two coordinates are the same, that is, $n = k$. This follow from (3), and therefore it is true for any projection and, in particular, for the Dyck *nk*-triangle. This means that in the Dyck *nk*-triangle the *j*-isoline #0 is the main diagonal. In this case, we get a lattice with $n \times n$ (or $k \times k$) square cells. Such square diagrams are well known [6, 7].

Let's transform the Dyck *ij*-triangle in Figure 2. We twirl the diagram around the main diagonal 180 degrees (change axes) and then rotate clockwise by 45 degrees. As a result, we formed the Dyck *nk*-triangle (see Figure 4). In a new projection, the diagonals became horizontal and vertical lines. The same Dyck *nk*-path is the monotonic lattice path along the edges of the *nk*-grid.

In the Dyck *nk*-triangle, the arbitrary Dyck *nk*-path starts in the lower left corner, finishes in the upper right corner, and consists entirely of edges pointing rightwards or upwards. Such a path from (0, 0) to $(n, n)$ (or $(k, k)$) does not rise above the main diagonal $n = k$.

In Figure 4, the usual square lattice is complemented by diagonal $i$-isolines and $j$-isolines for better orientation at the grid nodes. For example, it is easy to verify the dynamic equation (6) for three selected nodes.

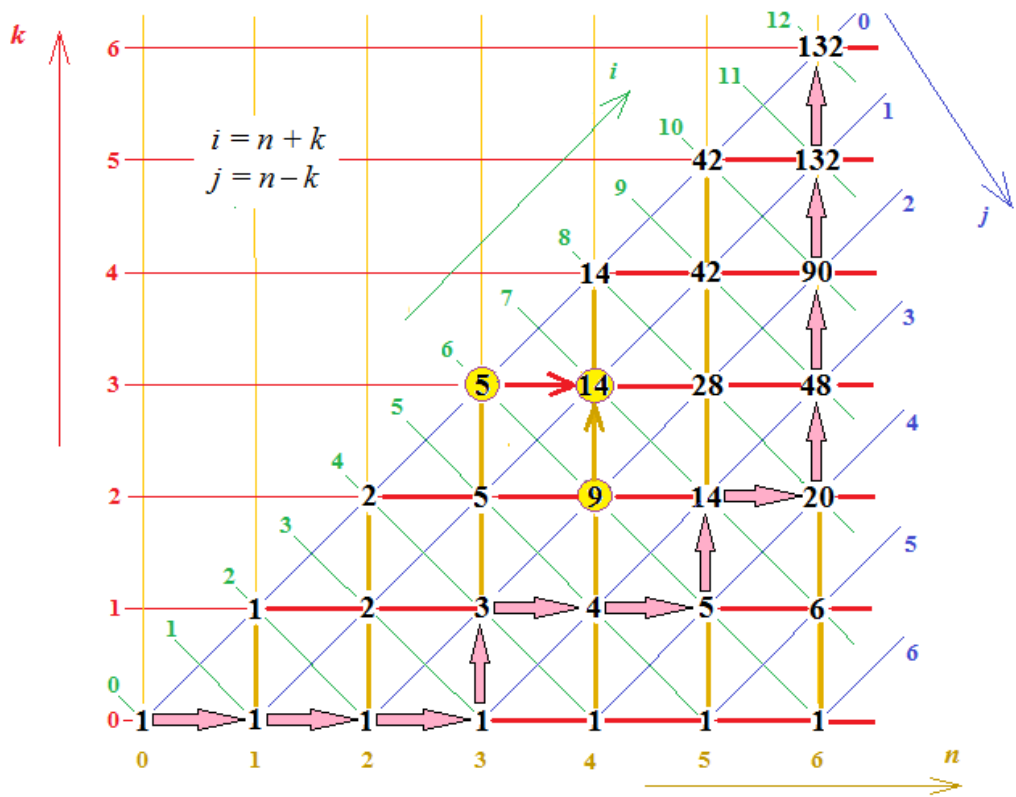

Figure 4: The Dyck $nk$-triangle.

Let's pay attention to the following. All three described triangles are well known and studied in detail. The first two grids $i \times j$ and $n \times j$ have the clear and natural coordinates. We have another situation with the classic square lattice $n \times n$. Here both coordinates are intuitively related to the Catalan number index, i.e. they do not seem to differ. In this case, on the abscissa we indicate the index of the associated Catalan number; the ordinate is the index of the Dyck square.

## 5 Other 2D Dyck triangles

In the Dyck $ijnk$-triangle, each node is uniquely determined by any pair of coordinate variables, and the other two coordinates are computed from (3). This means that all 2D and 3D projections are equivalent, mutually interchangeable. For example, we can easily move from the Catalan convolution matrix (the Dyck $nj$-triangle) with traced Dyck paths to the $ij$-triangle, on which we will see the same paths. Obviously, to solve a certain problems, this or that projection is preferred.

We have considered three known two-dimensional projections of the 4D Dyck triangle; other projections are not found in the literature. There are three more projections in the coordinate grids $i \times n$, $i \times k$, and $k \times j$ (in each grid, the order of coordinates is not important). It is logical to assume that these projections are not yet in demand.

**5.1. The Dyck *in*-triangle.** In Figure 5, we have traditionally indicated the position of the brackets along the abscissa axis, and the ordinate axis is the index of the associated Catalan number. The resulting Dyck *in*-triangle is more like a wedge, not a triangle. The *j*-isolines are not shown here to simplify the drawing. The lower nodes (Catalan numbers) of the *k*-isolines are on the *j*-isoline #0. The second group of Catalan numbers is on *j*-isoline #1.

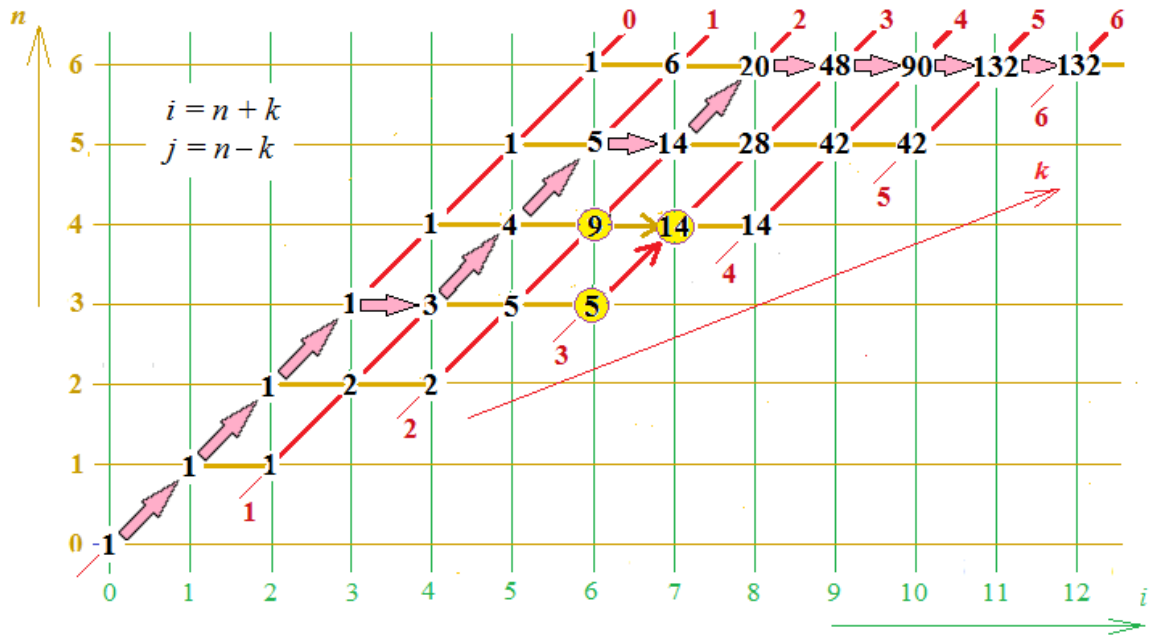

Figure 5: The Dyck *in*-triangle.

**5.2. The Dyck *kj*-triangle.** In Figure 6, we have traditionally indicated the unbalance of the brackets along the ordinate axis, and the abscissa axis is the index of the Dyck square. In fact, the resulting Dyck *kj*-triangle can no longer be called a triangle, since the first quadrant is fully accessible to Dyck paths (there are no forbidden areas).

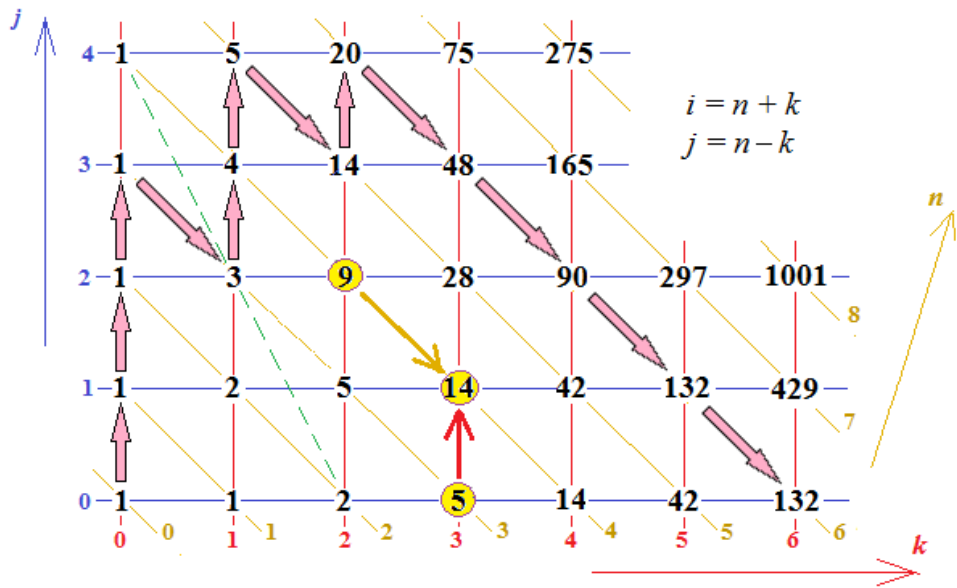

Figure 6: The Dyck *kj*-triangle.

It also does not show all isolines. Nodes that lie on the same *i*-isoline are located far from each other. For example, the green dashed line, *i*-isoline #4, links *kj*-nodes (0, 4), (1, 2), and (2, 0). Recall the respective coordinate equation $i = j + 2k$.

**5.3. The Dyck *ik*-triangle.** In Figure 7 on the abscissa axis, we indicated the usual bracket positions, and the ordinate axis is the indices of the Dyck square. The resulting picture is an unusual triangle. Not all isolines are shown in the *ik*-triangle. On the whole scheme we have drawn only the *j*-isoline #0 (dotted blue line) which connects the Catalan numbers. As we see, the Catalan numbers are the upper elements in the columns.

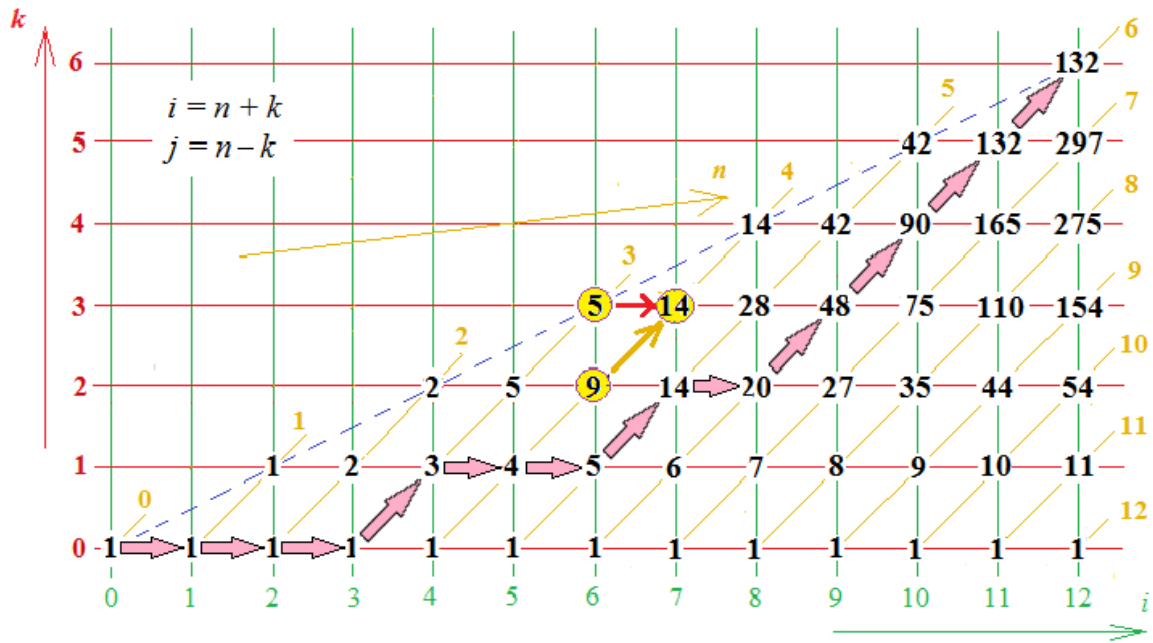

Figure 7: The Dyck *ik*-triangle

Let us pay attention once again, in all projections of the 4D Dyck triangle, 0-isolines intersect at the origin.

# 6 3D Dyck triangles

In this section, we construct four 3D projections of the 4D Dyck triangle. Obviously, in projections the orientation and order of the axes doesn't really matter. We will try to find a convenient image. In the general case, if possible, we will adhere to the order of the coordinates in the $i \times j \times n \times k$ lattice.

It is convenient to create a 3D projection based on some selected 2D triangle. For example, from the *ij*-triangle (Figure 2) we can get the *ijn*-triangle or the *ijk*-triangle. Accordingly, adding one coordinate in the *nk*-triangle (Figure 4), it is not difficult to construct the *ink*-triangle or the *jnk*-triangle.

For a better orientation in the three-dimensional projections, we will mark rays indicating all four coordinates of the points on these rays. This applies to the axes and to different diagonals. In the rays, certain coordinates of the points can be the same, some coordinates are fixed. In case of coincidence of coordinates, the first variable is repeated. Here are a few examples:

($i$, 0, 0, 0), (0, $j$, 0, 0), (0, 0, $n$, 0) и (0, 0, 0, $k$) – the main axes in the 4D lattice;
($i$, $i$, $i$, 0) – the central ray in the Dyck $ijn$-triangle (or $k$-diagonal #0);
($i$, $j$, $j$, 0) – the main diagonal in the Dyck $nj$-triangle (see Figure 3);
($i$, 4-$i$, 2, $i$-2) – the $n$-diagonal #2 in the Dyck $ijn$-triangle (see Figure 8).

In the drawings, the isolines are shown in different color:

$i$-isolines – green, $j$-isolines – blue, $n$-isolines – dark yellow, $k$-isolines – red.

In the projections, we will also repeat the coordinate equations (3), the old Dyck path, and the three nodes (6, 0, 3, 3), (6, 2, 4, 2), (7, 1, 4, 3). Let's repeat the dynamics equation for these nodes: $D(7, 1, 4, 3) = D(6, 0, 3, 3) + D(6, 2, 4, 2)$.

**6.1. The Dyck $ijn$-triangle.** The Figure 8 shows an $i \times j \times n$ lattice, that based on the Dyck $ij$-triangle. As we see the Dyck $ijn$-triangle is flat as all previous projections. This is not surprising since the coordinate equations (3) are linear.

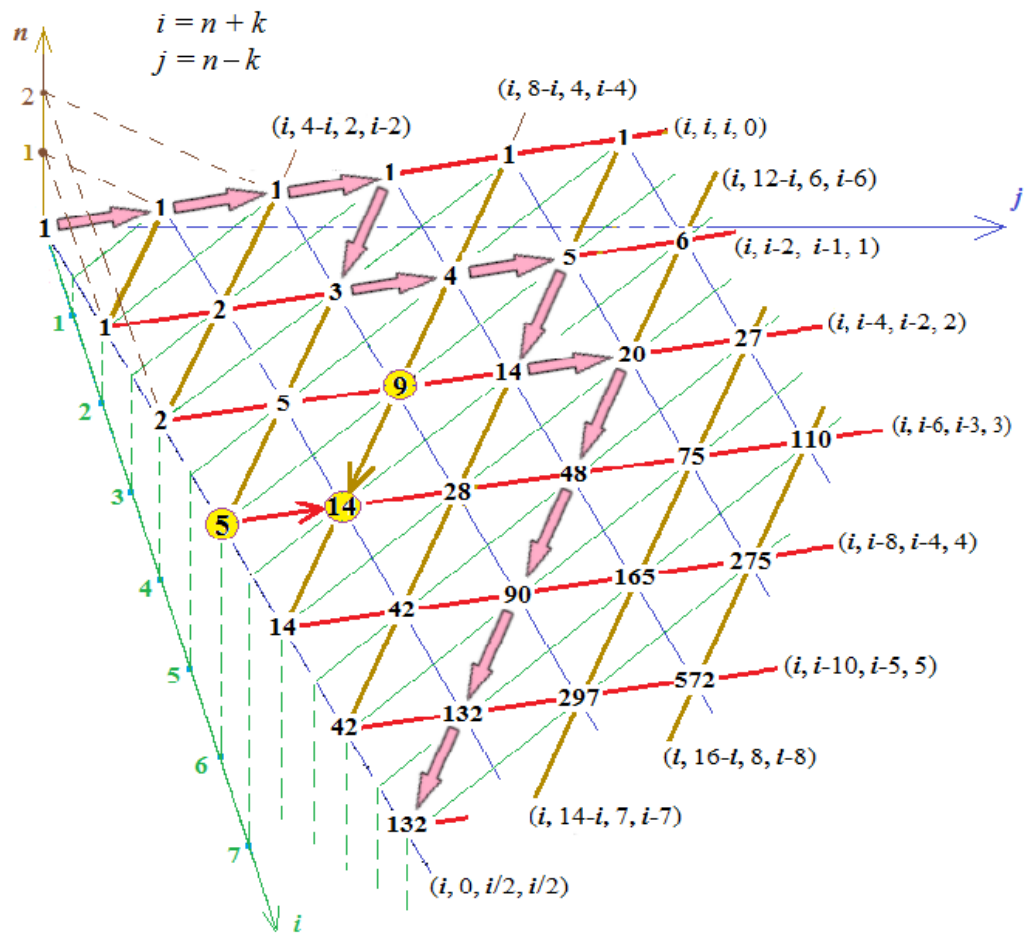

Figure 8: The Dyck $ijn$-triangle

**6.2. The Dyck $ijk$-triangle.** The Figure 9 shows an $i \times j \times k$ lattice, that also based on the Dyck $ij$-triangle. As we see the Dyck $ijk$-triangle is flat.

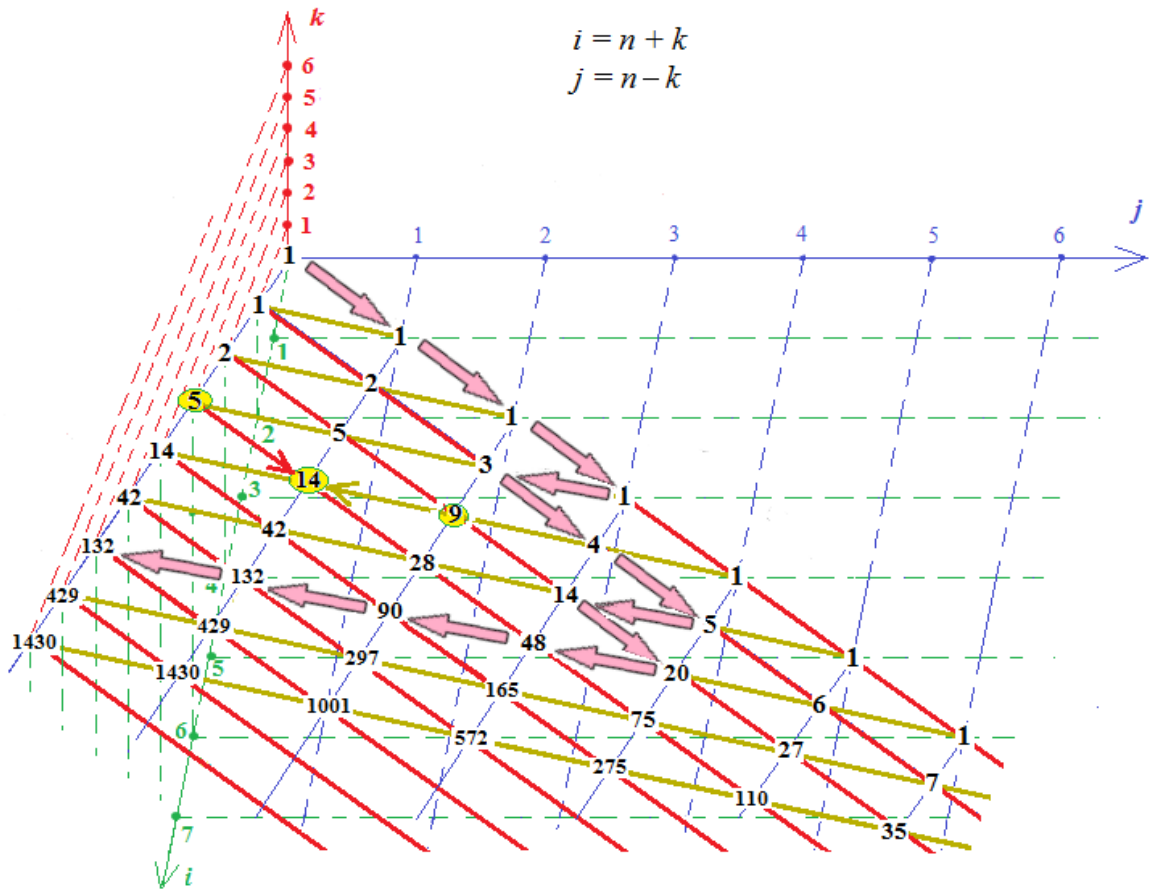

Figure 9: The Dyck *ijk*-triangle

**6.3. The Dyck *nik*-triangle.** The Figure 10 shows an $n \times i \times k$ lattice, that based on the Dyck *nk*-triangle (see Figure 4). This projection of the 4D Dyck triangle is also flat.

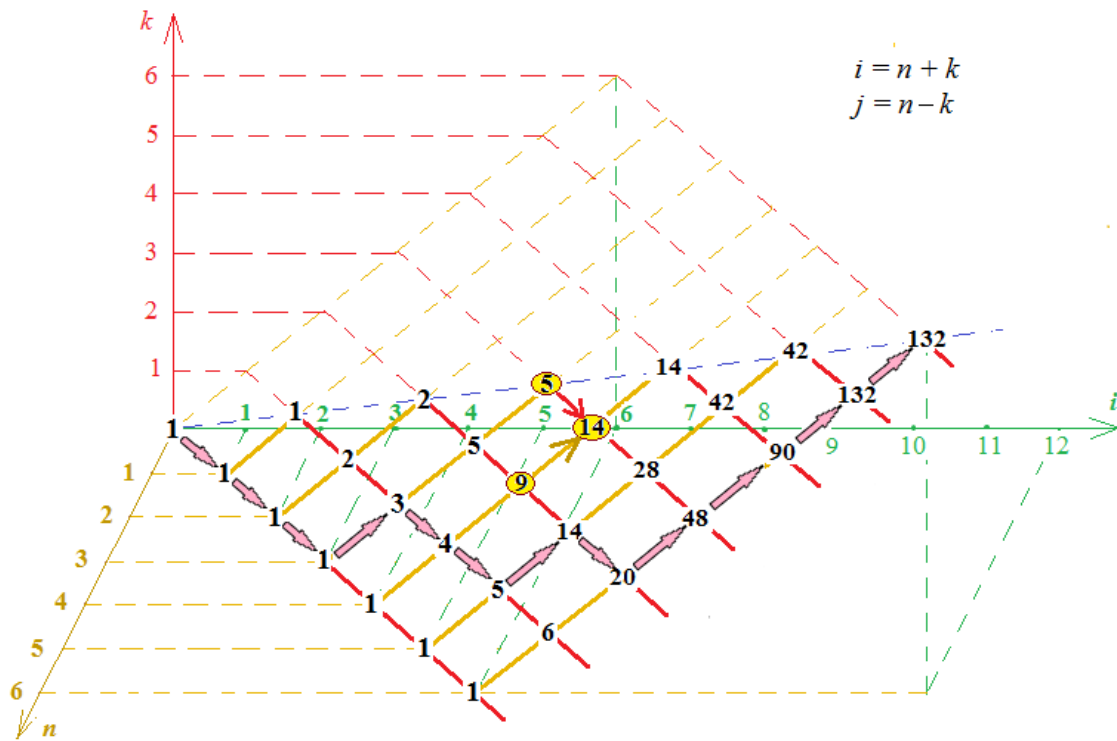

Figure 10: The Dyck *nik*-triangle

**6.4. The Dyck *jnk*-triangle.** The Figure 11 shows a $j \times n \times k$ lattice, that also based on the Dyck *nk*-triangle (see Figure 4). Figure 11 shows approximately the direction of the fourth i-axis. This projection of the Dyck *ijnk*-triangle is also flat.

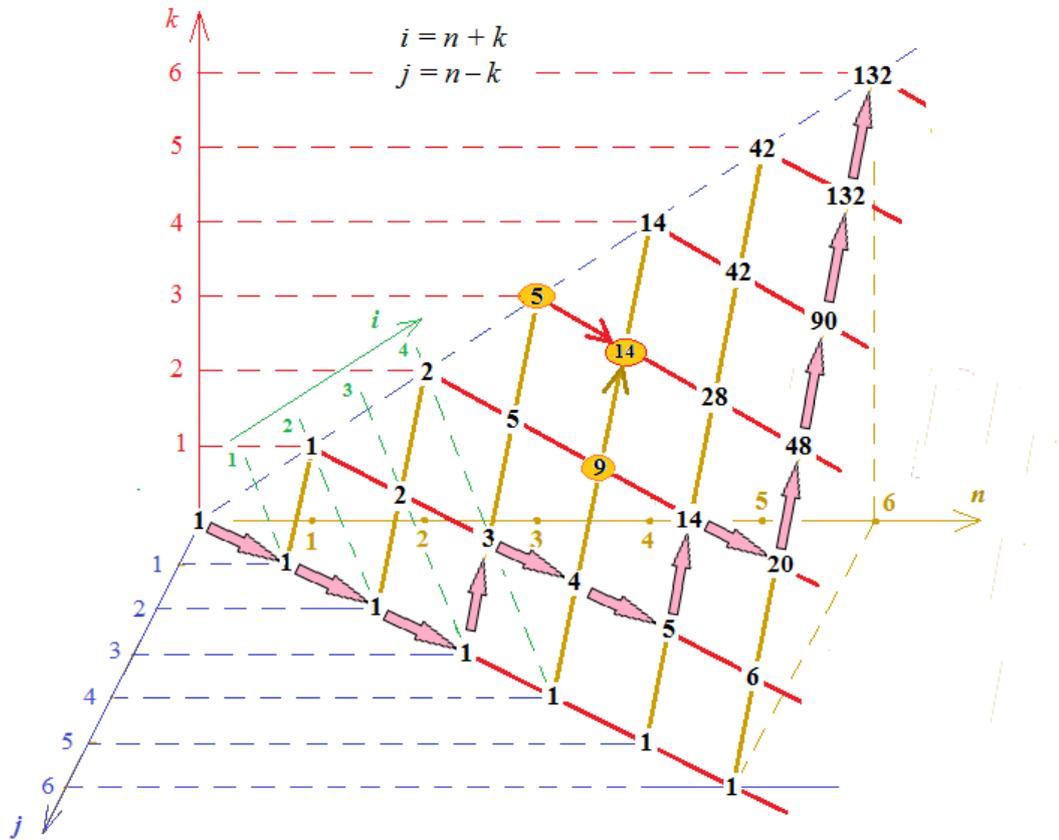

Figure 11: The Dyck *jnk*-triangle

Gzhel State University, Moscow, 140155, Russia
http://www.en.art-gzhel.ru/